\documentclass[a4j,12pt]{article}

\usepackage{amsmath}
\usepackage{amssymb}
\usepackage{amsthm}
\usepackage{graphics}
\usepackage[dvips]{graphicx}
\usepackage{amscd}
\usepackage[alphabetic,abbrev]{amsrefs}

\numberwithin{equation}{section}
\theoremstyle{definition}
\newtheorem{example}{example}[section]
\newtheorem{df}[example]{Definition}
\theoremstyle{remark}
\newtheorem{rem}[example]{Remark}
\theoremstyle{plain}
\newtheorem{lm}[example]{Lemma}
\newtheorem{prop}[example]{Proposition}

\newtheorem{thm}[example]{Theorem}

\theoremstyle{plain}
\newtheorem{theorem}[example]{Theorem}
\theoremstyle{remark}

\def\Real{{\Bbb R}}
\def\Complex{{\Bbb C}}

\def\Rc{\textrm{Ric}}

\theoremstyle{remark}

\begin{document}
\title{Ricci Flow on 3-dimensional Lie groups and 4-dimensional Ricci-flat manifolds}
% Ricci Flow on $3$dimension Lie Groups and
% $4$dimension Ricci-flat manifolds\\
\author{Kensuke Onda\footnote{kensuke.onda@math.nagoya-u.ac.jp}}
\date{Graduate School of Mathematics, Nagoya University }
\maketitle

% NOT of DISTRIBUTRION  
% Ryoichi Kobayashi,

\section*{Abstract}
We study relation of the Ricci Flow on $3$-dimensional Lie groups and $4$-dimensional Ricci-flat manifolds.
In particular, 
 we construct Ricci-flat cohomogeneity one metrics for $3$ dimensional Lie groups.
% group of flat Euclid plane and flat Lorentz plane.

\section{Introduction}
In this paper, we will discuss some relationship between the Ricci flow of left-invariant Riemannian metrics 
on 3-dimensional unimodular simply-connected Lie group $G$ 
and the Ricci-flat metrics of cohomogeneity one on the space-time $\Bbb R\times G$. 

First of all, we introduce the notion of cohomogeneity one metrics with respect to  a Lie group $G$. 

\begin{df}  
A pseudo-Riemannian manifold $(M, g)$ is a cohomogeneity one with respect to  a Lie group $G$, if and only if
$G$ is a subgroup of ${\rm Isom}(M, g)$, and the codimension of principal orbits under the action 
of $G$ equals $1$. 
\end{df}

The simplest example of a cohomogeneity one metric arises from the standard 
action of $SO(n)$ to $\Real^n$. 
The singular orbit is $\{ 0\}$, and the principal orbits are $S^{n-1}$ of various radii. 

In this paper, we attempt to construct cohomogeneity one Einstein metrics (in particular Ricci-flat metrics) 
from the Ricci flow solutions of left-invariant metrics on 3-dimensional unimodular simply connected Lie groups. 
The Ricci-flat metrics obtained in this paper has the property that their sectional curvatures decay to $0$ 
at one end toward which the metric is complete. 
Such metrics was studied by Lorentz, Gibbons, Hawking, Pope and many people. 
Some of them is called ALF and the Taub-NUT metric on $\Real\times SU(2)$ 
arising from the Ricci flow of the left $SU(2)$- and right $U(1)$-invariant metrics on $SU(2)$ 
is ALF. 
% The metric which has the condition that is stronger than ALF is ALE metric.
% The large class of ALE was discovered by Gibbons and Hawking.

Before proceeding to general 3-dimensional Lie groups, we examine the case of 
$SU(2)$, which is well known. 
There exists a left-invariant coframe $\{\theta^i\}_{i=1}^3$ on $SU(2)$ 
satisfying $d\theta ^i =2\theta ^j \wedge \theta ^k$, 
where $(i, j, k)$ are cyclic permutation of $\{1,2,3\}$. 
Then cohomogeneity one metrics with respect to  $SU(2)$ is described as
\begin{equation} \label{CohomogSU2}
g= dt^2+a (t)^2(\theta ^1)^2 +b (t)^2(\theta ^2)^2 +c (t)^2(\theta ^3)^2 .
\end{equation}
Some special cases of $(\ref{CohomogSU2})$ are listed in the following examples: 
\begin{enumerate}
\item
If $a = b = c$ are linear, then the metric $g$ has constant curvature $0$. 
\item
If $a$, $b$ and $c $ are constants, then the metric $g$ becomes a product metric on $\Real \times SU(2)$.
\item
If $a = b = c = \sin t$, then the metric $g$ has positive constant curvature. 
\item
If $a = b = c = \sinh t$, then the metric $g$ has negative constant curvature. 
\end{enumerate}
Thus, it is interesting to ask for which $\{a(t),b(t),c(t)\}$ the resulting cohomogeneity one metric 
is Einstein. 
In this paper, when the triple of functions $\{a(t),b(t),c(t)\}$ satisfies the Ricci flow equation, 
we examine whether the resulting cohomogeneity one metric is a Ricci-flat metric. 
This generalizes the above example 1. 

First of all, we define the {\it Ricci flow} again. 
In this paper, we define that 
a 1-parameter family $g(t)$ of Riemannian metrics is the Ricci flow if and only if it solves
$$\frac{\partial }{\partial t} g(t)_{ij} = -\Rc [g (t)]_{ij} ,$$
and a 1-parameter family of Riemannian metrics $g(t)$ is the backward Ricci flow if and only if it solves
$$\frac{\partial }{\partial t} g(t)_{ij} = \Rc [g (t)]_{ij} .$$
In the Introduction, the Ricci flow was defined by $$\frac{\partial }{\partial t} g(t)_{ij} = -2\Rc [g]_{ij} .$$
The difference in $1$ and $2$ is just a choice of the scale and produces no trouble. 
Indeed, if we put $h: = \frac{k}{2} g$, where $k$ is a positive constant, then 
$$\frac{\partial }{\partial t} h(t)_{ij} = (\frac{k}{2} ) (-2\Rc [g]_{ij}) = -k \Rc [h]_{ij} , $$
because $\Rc [h]=\Rc [\frac{k}{2} g] = \Rc[g]$.
In this paper, if $k$ equals $1$, we say that $h(t)$ is a solution to the Ricci flow. 
Also the backward Ricci flow (in this paper) is defined similarly. 
Note that if we change $t$ into $-t$, the Ricci flow equation changes 
into the backward Ricci flow equation.  

Next, we review the Ricci flow on $SU(2)$.
Let $\{ F_i\} _{i=1}^3$ be an orthonormal frame of a left-invariant metric $g_3$ on $SU(2)$, 
satisfying $[F_i, F_j]=-2 F_k$  \ $(i, j, k:$ cyclic$)$, 
and $\{ \theta ^i\}$ the dual coframe of $\{ F_i\} _{i=1}^3$.
The left-invariant metric $g_3$ is expressed as
$$
g_3 = A (\theta ^1)^2+ B (\theta ^2)^2 +C (\theta ^3)^2 . 
$$
Then the Ricci flow is equivalent to the system of ODE's: 

\begin{equation} \label{SU2}
    \begin{aligned}
\frac{d}{dt} A = \frac{(B-C)^2-A^2}{B C} , \\
\frac{d}{dt} B=\frac{(C-A)^2-B^2}{C A} , \\  
\frac{d}{dt} C=\frac{(A-B)^2-C^2}{A B} .
     \end{aligned}
\end{equation}
The behavior of solutions of $(\ref{SU2})$ is known in \cite{KM01, CK}. 
\begin{prop}[\cite{KM01, CK}]
{\it The solution of the Ricci flow equation exists on $(-\infty, T)$, where $T$ depends on the initial data, 
and if $t$ goes to $T$, then $g_3$ becomes asymptotically round and shrinks to a point.}
\end{prop}
Next we consider a cohomogeneity one metric with respect to  $SU(2)$, given by 
\begin{equation} \label{CSU2}
g= dt^2+a (t)^2 (\theta ^1)^2 +b (t)^2 (\theta ^2)^2 +c (t)^2 (\theta ^3)^2 .
\end{equation}
A typical example of cohomogeneity one metric with respect to  
$SU(2)$ is the Taub-NUT metric, given by 
$$g = \Big( \frac{r+m}{r-m} \Big)dr^2 +(r^2-m^2)\{ (\theta^1)^2 +(\theta ^2)^2\} 
+4 m^2\Big( \frac{r-m}{r+m} \Big) (\theta ^3)^2  .$$
As is well known, this metric has the following properties. 
First, this metric is a hyper-K\"{a}hler metric, and therefore this is a Ricci-flat metric.
Secondly, we can put 
$$a = b = (r^2-m^2)^{\frac{1}{2}} , \quad c = 2 m\Big( \frac{r-m}{r+m} \Big)^{\frac{1}{2}}  , $$
and therefore the coefficient $a$ equals $b$, but $a$ is not equal to $c$ (i.e., the metric on 
$SU(2)$ part is left $SU(2)$- and right $U(1)$-invariant). 
Thirdly, the change of variables 
$h = \Big( \frac{r+m}{r-m} \Big)^{\frac{1}{2}}$ and $dt = h dr$ implies that 
this metric is regarded as a cohomogeneity one metric. 
Fourthly, it is easy to check that $a$, $b$ and $c$ satisfy 
\begin{equation*}
    \begin{aligned}
\frac{d}{dt} a & = \frac{a^2-(b-c)^2}{b c} , \\
\frac{d}{dt} b & =\frac{b^2-(c-a)^2}{c a} , \\
\frac{d}{dt} c & =\frac{c^2-(a-b)^2}{a b} .
     \end{aligned}
\end{equation*}
This system is equivalent to the backward Ricci flow equation for left-invariant metrics on $SU(2)$. 
Thus we conclude that if coefficients $\{a(t),b(t),c(t)\}$ move along the Ricci flow or the backward Ricci flow 
on $SU(2)$, then the resulting cohomogeneity one metric $(\ref{CSU2})$ is the Taub-NUT metric and therefore 
Ricci-flat.

\begin{thm}[\cite{F85, CGLP04}]\label{su2ricciflat} %G. W. Gibbons, '98]
If $a$, $b$ and $c$ satisfy the Ricci flow equations or the backward Ricci flow equations of a metric 
\begin{equation}\label{3metric}
g_3=a(\theta ^1)^2+b(\theta ^2)^2+c(\theta ^3)^2
\end{equation}
on $SU(2)$, 
then the cohomogeneity one metric
\begin{equation}\label{cohomo_su2}
g= dt^2+a(t)^2(\theta ^1)^2 +b(t)^2(\theta ^2)^2 +c(t)^2(\theta ^3)^2
\end{equation}
with respect to  $SU(2)$ on the space-time of the Ricci flow becomes a Ricci-flat metric. 
\end{thm}
It is important to pay attention to coefficients. 
Coefficients of metric $(\ref{3metric})$ are $a$, $b$ and $c$.
But, coefficients of metric $(\ref{cohomo_su2})$ are $a^2$, $b^2$ and $c^2$.
Therefore, $a$, $b$ and $c$ are positive (because we consider left-invariant Riemannian metrics 
on $SU(2)$). Even if $a$, $b$ and $c$ are positive, we have still freedom in introducing minus sign 
before $\{a^2,b^2,c^2\}$ in the attempt to discover cohomogeneity one metrics on the space-time of 
the Ricci flow (this freedom is really essential in the case of $E(1,1)$ and $SL(2,\Real)$). 

\begin{rem}
Even if a cohomogeneity one metric
$$
g= dt^2+ a (t)^2(\theta ^1)^2 +b (t)^2(\theta ^2)^2 + c (t)^2(\theta ^3)^2
$$
with respect to  $SU(2)$ is Ricci-flat, it may not satisfy the Ricci flow equation 
nor the backward Ricci flow equation. 
We consider the Eguchi-Hanson metric as a typical example.
The Eguchi-Hanson metric is given by
\begin{equation}\label{eguchi}
g=\frac{dr^2}{1-(m/r)^4} +r^2\{ (\theta^1)^2 +(\theta ^2)^2\} +r^2( 1-(m/r)^4) (\theta ^3)^2.
\end{equation}
Put 
$$a=b=r , \ c=r\Big( 1-(\frac{m}{r} )^4\Big)^{\frac{1}{2}}  .$$
By the coordinate transform, 
$$h=\frac{1}{\Big( 1-(\frac{m}{r} )^4\Big)^{\frac{1}{2}}} , \  dt=h dr , $$
the metric $(\ref{eguchi})$ is regarded as a cohomogeneity one metric.
Also this metric $(\ref{eguchi})$ is a Ricci-flat metric,
but coefficients $a$, $b$ and $c$ satisfy the following system of ODE's: 
\begin{equation*}
    \begin{aligned}
\frac{d}{d t} a & = \frac{a^2-( b - c )^2}{b c} +2, \\
\frac{d}{d t} b & = \frac{b^2-( c - a )^2}{c a} +2, \\
\frac{d}{d t} c & = \frac{c^2-( a - b )^2}{a b} +2.
	     \end{aligned}
\end{equation*}
These equations are neither the Ricci flow nor the backward Ricci flow. 
\end{rem}

We now study cohomogeneity one metrics constructed from the Ricci flow solution on other groups, 
and construct Ricci-flat metrics on their space-time. 
Ricci-flat metrics that we construct in this paper are listed in Table 1.

\begin{table}[htbp]
\begin{center}
\begin{tabular}{|c|c|c|}
\hline
Lie groups & the Ricci flow equation & signature of metric \\ \hline
$SU(2)$ & $SU(2)$ & $(4, 0)$ \\ \hline
$E(2)$ & $E(2)$ & $(4, 0)$ \\ \hline
$SL(2, \Real )$ & $SU(2)$ & $(2, 2)$ \\ \hline
$H_3$ & $H_3$ & $(4, 0)$ \\ \hline
$E(1, 1)$ & $E(1, 1)$ & $(2, 2)$ \\ 
\hline
\end{tabular}
\caption{Ricci-flat metrics of cohomogeneity one corresponding to 3-dimensional Lie groups}
\end{center}
\end{table} 

% This paper is organized as follows.  %teisei 
% In Subsection 1, 
% we introduce a Milnor frame and calculate the Ricci tensor of cohomogeneity one metric to prepare for 
% Subsection 3-6.
% In Subsection 2, we construct a cohomogeneity one Ricci-flat metric on the space-time of the Ricci flow 
% of left-invariant Riemannian metrics on the Heisenberg group $H_3$ 
% and see that this metric has the hyper-K\"{a}hler structure.
% In Subsection 3,  
% we construct a sign $(2, 2)$ cohomogeneity one Ricci-flat metric on the space-time of the Ricci flow 
% of left invariant Riemannian metrics on the group $E(1,1)$ of rigid motions of the Minkowski $2$-space. 
% In Subsection 4, 
% we construct a cohomogeneity one Ricci-flat metric on the space-time of the Ricci flow 
% of left-invariant Riemannian metrics on the group of rigid motions of the Euclidean $2$-space $E(2)$ 
% and see that this metric has the hyper-K\"{a}hler structure. 
% In Subsection 5, we construct a sign $(2, 2)$ cohomogeneity one Ricci-flat metric on the space-time 
% of the Ricci flow of left-invariant Riemannian metrics on $SL(2, \Real )$.

\section{Preparation}
We present the definition of a Milnor frame.
\begin{df} 
Let $\{ F_i\} _{i=1}^3$ be a left-invariant moving frame on $G$.
If $\{ F_i\} _{i=1}^3$ satisfies 
$$[F_2, F_3] = n_1 F_1, \quad [F_3, F_1] = n_2 F_2, \quad [F_1, F_2] = n_3 F_3 , $$
where $\ n_i \in \{ \pm 1, 0\} $, 
then $\{ F_i\} _{i=1}^3$ is called {\it a Milnor frame} .
\end{df}
As is well known, 3-dimensional unimodular simply-connected Lie groups were classified by Milnor \cite{M76}.
\begin{prop}[\cite{M76}]
{\it Let $\{ F_i\}$ be a Milnor frame. 
For signatures of $\{ n_i\}$, $3$-dimensional unimodular simply-connected Lie groups are determined as Table $2$.}
\begin{table}[htbp]
\begin{center}
\begin{tabular}{|c|c|c|}
\hline
Signature & Lie groups & description \\ \hline
$(-1, -1, -1)$, $(+1, +1, +1)$ & $SU(2)$ & simple \\ \hline
$(-1, -1, +1)$, $(-1, +1, +1)$ & $\widetilde{{\rm SL}(2, \Real )}$ & simple \\ \hline
$(-1, -1, 0)$, $(+1, +1, 0)$ & $E(2)$ & solvable \\ \hline
$(-1, 0, +1)$ & $E(1, 1)$ & solvable \\ \hline
$(-1, 0, 0)$, $(+1, 0, 0)$ & $H_3$ & nilpotent \\ \hline
$(0, 0, 0)$ & $\Real \oplus \Real \oplus \Real$ & commutative \\ \hline
\end{tabular}
\caption{$3$-dimensional unimodular simply-connected Lie groups}
\end{center}
\end{table} 
\end{prop}

\begin{rem}
We may change order of $n_i$.
\end{rem}

Let $(M^4, g)$ be a cohomogeneity one 4-dimensional manifold with respect to $3$-dimensional Lie group $G$. 
Let $\{ F_i\} _{i=1}^3$ be a Milnor frame of $G$, and 
$\{ \theta _i\}$ the dual coframe of $\{ F_i\}$. 
Then the metric is expressed as follows:
\begin{equation}\label{cohomo1}
g=dt^2+a(t)^2 (\theta ^1)^2+b(t)^2 (\theta ^2)^2+c(t)^2 (\theta ^3)^2 .
\end{equation}
From Theorem $1.1$ (Koszul's formula), the Levi-Civita connection is expressed as below.
\begin{prop} \label{levi2}
{\it Let $g$ be the cohomogeneity one metric $(\ref{cohomo1})$. Then the Levi-Civita connection is given by}
\begin{eqnarray}
% (\nabla _{F_i} F_j)=
\left(
\begin{array}{cccc}
0 & \frac{\dot{a}}{a} F_1 & \frac{\dot{b}}{b} F_2 & \frac{\dot{c}}{c} F_3 \\
\frac{\dot{a}}{a} F_1 & -a\dot{a} F_0 & \frac{1}{2} \cdot \frac{-n_1 a^2+n_2 b^2+n_3 c^2}{c^2} F_3 & \frac{1}{2} \cdot \frac{-n_3 c^2+n_1 a^2-n_2 b^2}{b^2} F_2 \\
\frac{\dot{b}}{b} F_2 & \frac{1}{2} \cdot \frac{-n_1 a^2+n_2 b^2-n_3 c^2}{c^2} F_3 & -b\dot{b} F_0 & \frac{1}{2} \cdot \frac{-n_2 b^2+n_3 c^2+n_1 a^2}{a^2} F_1 \\
\frac{\dot{c}}{c} F_3 & \frac{1}{2} \cdot \frac{-n_3 c^2+n_1 a^2+n_2 b^2}{b^2} F_2 & \frac{1}{2} \cdot \frac{-n_2 b^2+n_3 c^2-n_1 a^2}{a^2} F_1 & -c\dot{c} F_0  \\
\end{array}
\right) . 
\end{eqnarray} 
\end{prop}
Hence we obtain the Ricci tensor; 
\begin{subequations}\label{cohomo_rc}
    \begin{align}
\Rc (F_0, F_0)=R_{00} & = -\frac{\ddot{a}}{a} -\frac{\ddot{b}}{b} -\frac{\ddot{c}}{c}  , \\
\Rc (F_1, F_1)=R_{11} & = -a \frac{(\dot{a} b c)^{\cdot }}{b c} -\frac{(n_2 b^2-n_3 c^2)^2-n_1^2 a^4}{2 b^2 c^2} , \\
\Rc (F_2, F_2)=R_{22} & = -b \frac{(\dot{b} c a)^{\cdot }}{c a} -\frac{(n_3 c^2-n_1 a^2)^2-n_2^2 b^4}{2 c^2 a^2} , \\
\Rc (F_3, F_3)=R_{33} & = -c \frac{(\dot{c} a b)^{\cdot }}{a b} -\frac{(n_1 a^2-n_2 b^2)^2-n_3^2 c^4}{2 a^2 b^2} , 
 \end{align}
 \end{subequations}
and other components are $0$.
Remark that if $a$, $b$ and $c$ are constants, then the Ricci tensor $R_{11}$, $R_{22}$ and $R_{33}$ of the metric $(\ref{cohomo1})$ are the Ricci tensors of a metric
$$g=a^2 (\theta ^1)^2+b^2 (\theta ^2)^2+c^2 (\theta ^3)^2 .$$ 

%%%%%%%%%NILPOTENT-GEOMETRY, nil, Nil, Heisenberg
\section{The Heisenberg group}
In this section, we consider the Heisenberg group $H_3$ and a cohomogeneity one with respect to $H_3$.
Let $\{ F_i\}$ be the Milnor frame of $H_3$, 
satisfying $[F_1, F_2] = 0$, $[F_2, F_3] = F_1$ and $[F_3, F_1] = 0$; 
and $\{ \theta _i\}$ the dual coframe of $\{ F_i\}$.
Then a left-invariant metric is expressed as
\begin{equation}\label{nil3}
g = a (\theta ^1)^2 + b (\theta ^2)^2 + c (\theta ^3)^2 .  
\end{equation}
Therefore the Ricci flow equation is given by
\begin{subequations}\label{nilrf}
    \begin{align}
    \frac{d}{dt} a = & -\frac{a^2}{2 b c} , \\
    \frac{d}{dt} b = & \frac{a}{2 c} , \\
    \frac{d}{dt} c = & \frac{a}{2 b} .
   \end{align}
\end{subequations}
These equations were solved in \cite{KM01} as follows: 
\begin{equation}\label{nilsolu}
    \begin{cases}
     a (t) = a_0^{\frac{2}{3}} b_0^{\frac{1}{3}} c_0^{\frac{1}{3}} (\frac{3}{2} t + b_0 c_0/a_0 )^{-\frac{1}{3}} , \\
    b (t) = a_0^{\frac{1}{3}} b_0^{\frac{2}{3}} c_0^{-\frac{1}{3}} (\frac{3}{2} t + b_0 c_0/a_0 )^{\frac{1}{3}} ,  \\
    c (t) = a_0^{\frac{1}{3}} b_0^{-\frac{1}{3}} c_0^{\frac{2}{3}} (\frac{3}{2} t + b_0 c_0 /a_0 )^{\frac{1}{3}} ,    
     \end{cases}
\end{equation}
where $a_0=a(0)$, $b_0=b(0)$ and $c_0=c(0)$.  
In particular, the behavior of this solution is the following: 

\begin{lm}[\cite{KM01}]\label{heisenRF}
{\it The metric $(\ref{nil3})$ satisfying $(\ref{nilrf})$ has three properties.}
\begin{enumerate}
\item {\it The solution of the Ricci flow equation exists on} $(-T,  +\infty)$, {\it where} $T=\frac{b_0 c_0}{3 a_0} $.
\item {\it If} $t \rightarrow -T$, 
 {\it then $a \rightarrow +\infty$, $b \rightarrow 0$ and $c \rightarrow 0$. }
\item {\it If $t \rightarrow +\infty$, 
then $a \rightarrow 0$, $b \rightarrow +\infty$ and $c \rightarrow +\infty$.}
\end{enumerate}
% \end{itemize}
\end{lm}
We use this lemma later to describe the resulting Ricci-flat metric on the space-time. 
Next we consider a cohomogeneity one metric with respect to  the Heisenberg group $H_3$. 
Let $(M^4, g) $ be a cohomogeneity one with respect to $H_3$. 
% Let $\{ F_i\} _{i=1}^3$ and $\{ \theta ^i\}$ be as beofe.
% Then metric is expressed as 
Using the above frame, a cohomogeneity one metric is described as
\begin{equation} \label{cohomo_hisen}
g=dt^2 + a (t)^2 (\theta ^1)^2 + b (t)^2 (\theta ^2)^2 + c (t)^2 (\theta ^3)^2 .
\end{equation}
Put $F_0 =\frac{\partial}{\partial t}$, then $\{ F_i\}$ becomes an orthogonal frame.
The Ricci tensor of this metric is given by
\begin{subequations}
    \begin{align*}
    R_{00} & = -\frac{\ddot{a}}{a} -\frac{\ddot{b}}{b} -\frac{\ddot{c}}{c} , \\
    R_{11} & = -a \frac{(\dot{a} b c)^{\cdot }}{b c} + \frac{a^4}{2 b^2 c^2} , \\
    R_{22} & = -b \frac{(\dot{b} c a)^{\cdot }}{c a} -\frac{a^2}{2 c^2} , \\
    R_{33} & = -c \frac{(\dot{c} a b)^{\cdot }}{a b} -\frac{a^2}{2 b^2} , 
 \end{align*}
 \end{subequations}
and other components are $0$.
Using $(\ref{nilsolu})$, we construct a Ricci-flat metric in the following way: 
\begin{thm} [\cite{GR98}]
{\it Let $a$, $b$ and $c$ satisfy the Ricci flow equations on $H_3$: 
\begin{subequations}
  \begin{align*}
    \frac{d}{dt} a & =-\frac{a^2}{2 b c}, \\
    \frac{d}{dt} b & =\frac{a}{2 c}, \\
    \frac{d}{dt} c & =\frac{a}{2 b}.
    \end{align*}
 \end{subequations}
Then $g$ becomes a Ricci-flat metric.} 
\end{thm}

\begin{proof}
An easy computation shows that 
\begin{subequations}
  \begin{align*}
    \frac{d^2}{dt^2} a & = \frac{a^3}{b^2 c^2} , \\
    \frac{d^2}{dt^2} b & =- \frac{a^2}{2 b c^2} , \\
    \frac{d^2}{dt^2} c & =- \frac{a^2}{2 b^2 c} , 
    \end{align*}
 \end{subequations}
and $a b$, $a c$ and $b/c$ are constants. 
Using this fact and the equation $(\ref{cohomo_rc})$, we obtain the statement.
\end{proof}

Sectional curvatures of this metric are
$$
K_{01} = K_{23} = -\frac{a^2}{b^2 c^2} , \quad K_{02} =K_{13} = \frac{a^2}{2 b^2 c^2} , \quad  
K_{03} =K_{12} = \frac{a^2}{2 b^2 c^2} .  
$$
From Proposition $(\ref{heisenRF})$, we get the following proposition.
\begin{prop}
{\it Assume that $a_0 b_0 c_0$ is positive. Then}
% \begin{itemize}
\begin{enumerate}
\item {\it the parameter $t$ exists on $(-T, \infty)$, 
where} $T = \dfrac{2 b_0 c_0}{ 3 a_0} $.
\item {\it If} $t \rightarrow -T $, 
{\it then} $K_{ij} \rightarrow \pm \infty$.
\item {\it If} $t \rightarrow +\infty$, {\it then} $K_{ij} \rightarrow 0$.
\end{enumerate}
{\it Assume that} $a_0 b_0 c_0$ {\it is negative. Then} 
\begin{enumerate}
\item {\it the parameter $t$ exists on $(-\infty, T' )$, 
where} $T' = - \dfrac{2 b_0 c_0 }{3 a_0} $.
\item If $t \rightarrow -\infty $, {\it then} $K_{ij} \rightarrow T'$.
\item {\it If} $t \rightarrow +\infty$, {\it then} $K_{ij} \rightarrow 0$.
\end{enumerate}
\end{prop}
% \begin{rem}
% cohomogeneity one of Heisenberg group has $2$ parameter family in next sense.
% We consider $ac=\lambda $ and $b$.
% If $b \rightarrow \infty$, then scalar curvature converges $0$.
% Thereofe Ricci flow solution exists $(-\infty, \infty )$, and $(M^4, g)$ becomes complete, because it is  near Einstein.
% \end{rem}
Actually, this phenomenon is related with the hyper-K\"{a}hler structure.
We confirm it from now on.
Before we define almost complex structures, 
we put 
$$e^0 = d t, \quad e^1=a\theta ^1, \quad e^2 = b\theta ^2, \quad e^3 = c\theta ^3 .$$
Then $\{ e^i\} _{i=0}^3$ becomes an orthonormal coframe.
Let $\{ e_i\}$ be the orthonormal frame defined by $ e^i(e_j)=\delta ^i_j$.
For $i=1$, $2$ and $3$, we consider almost complex structures $J_i$ defined by
$$J_i e_0 = e_i , \quad J_i e_j = - e_k , \quad (J_i)^2 = -id,$$ 
where $(i, j, k)$ are cyclic permutation of $\{1,2,3\}$. 
It is easy to check that $\{J_i\}$ satisfy $J_3 = -J_1 J_2 = J_2 J_1 .$ 
Moreover from direct calculation we have $N_{J_i}(F_j,F_k)=0$ for all $i,j,k=0,1,2,3$ 
where $N_{J_i}$ is the Nijenhuis tensor of the almost complex structure $J_i$. 
%%%%%%%%%%%%%%%%%%%%%%%%%%%%%%%%%%%%%%%%%%%%%%%%%%%%%%%
% You must show that the Nijenhuis tensor vanishes. 
% The Nijenhuis tensor of an almost complex structure $J$ is 
% a (1,2)-teisor dedfined by
% $$
% N(X,Y):=[JX,JY]-[X,Y]-J[JX,Y]-J[X,JY]
% $$
% and $N$ vanishes if and only if $J$ is integrable (i.e., holomorphic coordinate system exists locally). 
% You must compute $N(e_i,e_j)$ or $N(F_i,F_j)$ and show that these are zero. These are simple 
% computation. 
% Same remark in the $E(2)$-case also !! ok
%%%%%%%%%%%%%%%%%%%%%%%%%%%%%%%%%%%%%%%%%%%%%%%%%%%%%%%%
Therefore the almost complex structures $\{J_i\}$ are all integrable. 
We consider a triple of 2-forms $(\omega ^1, \omega ^2, \omega ^3)$, defined by 
$$\omega ^i (X, Y) := g (J_i X, Y)$$
for $i=1, 2, 3$.
Then 2-forms $\omega _i$ are given by
$$\omega ^i = e^0\wedge e^i - e^j \wedge e^k,$$
where $(i, j, k)$ are cyclic permutation of $\{1,2,3\}$. 
The above triple of $2$-forms $(\omega ^1, \omega ^2, \omega ^3)$ satisfies the following relations
$$(\omega ^i)^2\ne 0, \quad \omega ^i \wedge \omega ^j =0 \quad (i \ne j) .$$
Since $d e^0 = d (d t )=0, \quad  d e^1 = \dot{a} d t \wedge \theta ^1 - a \theta ^2 \wedge \theta ^3,
 \ d e^2 = \dot{b} d t \wedge \theta ^2, \quad d e^3 = \dot{c} d t \wedge \theta ^3,  $
 we get 
 \begin{subequations}
    \begin{align*}
d\omega ^1 & =(a -\dot{b} c -b\dot{c} ) d t \wedge \theta ^2 \wedge \theta ^3,  \\
d\omega ^2 & =(\dot{a} c +a \dot{c} ) d t \wedge \theta ^3 \wedge \theta ^1,  \\
d\omega ^3 & =(\dot{a} b +a \dot{b} ) d t \wedge \theta ^1 \wedge \theta ^2.
\end{align*}
\end{subequations}

\begin{thm}
Assume that the cohomogeneity one metric $(\ref{cohomo_hisen})$ has almost complex structures $J_i$. 
For any $i=1, 2, 3$,  $d \omega ^i =0$ if and only if
$a$, $b$ and $c$ satisfy the Ricci flow equations on the Heisenberg group. 
\end{thm}

\proof
We see at once that $d \omega ^i =0$ if and only if 
\begin{subequations}
    \begin{align}
a - \dot{b} c - b \dot{c} = 0 \label{na} ,  \\
\dot{a} c + a \dot{c} = 0 \label{nb} , \\
\dot{a} b + a \dot{b} = 0 \label{nc} .
\end{align}
\end{subequations}
Equations $(\ref{nb} )$ and $(\ref{nc} )$ imply that $a c$ and $a b$ are constants.
Using equations $(\ref{na} )$-$(\ref{nc})$, we compute that 
\begin{subequations}
    \begin{align*}
a = (b c)^{\cdot }  = & \big( \frac{a b\cdot a c}{a^2}\big) ^{.} \\ 
 = & \big( \frac{1}{a^2} \big) ^{\cdot } a^2 b c \\
  = & -2\dot{a} \frac{b c}{a} .
\end{align*}
\end{subequations}
Hence we get 
 $$\dot{a} =-\frac{a^2}{2 b c} .$$
Since
$ a\dot{b} =-b\dot{a} =a^2/(2 c)$ and $a\dot{c} =-c\dot{a} =a^2/(2 b) ,$   
we get 
\begin{subequations}
    \begin{align}
    \dot{b} = & \frac{a}{2 c} , \\
    \dot{c} = & \frac{a}{2 b} .
    \end{align}
\end{subequations}
Therefore $a$, $b$ and $c$ satisfy the Ricci flow equations on $H_3$. 

Conversely, $a, b$ and $c$ satisfy the Ricci flow equations on $H_3$, then $d \omega _i=0$.
% because propositon
\endproof
Therefore the cohomogeneity one metric with $(\ref{nilrf})$ becomes the hyper-K\"{a}hler metric.

%%%%%%%%%%%%%%%%%%%%%%%%%%Sol-Geometry, sol, Sol,

\section{The group of rigid motions of the Minkowski $2$-space}
In this section, we consider the group of rigid motions of the Minkowski $2$-space.
Let $\{ F_i\} _{i=1}^3$ be the Milnor frame of $E(1, 1)$, 
satisfying $[F_1, F_2] = - F_3, \ [F_3, F_1] = 0 , \ [F_2, F_3] = F_1$, 
and $\{ \theta _i\} $ the dual coframe of $\{ F_i\}$. 
Then a left-invariant metric is expressed as 
\begin{equation}\label{sol}
g = a (\theta ^1)^2 + b (\theta ^2)^2 + c (\theta ^3)^2 .
\end{equation}
Therefore the Ricci flow equation is given by
\begin{equation}\label{sol3rf}
    \begin{aligned}
    \dfrac{d}{dt} a & =\dfrac{c^2-a^2}{2 b c} , \\
    \dfrac{d}{dt} b & =\dfrac{( a + c )^2}{2 a c} , \\
    \dfrac{d}{dt} c & =\dfrac{ a^2 - c^2}{2 a b} . 
     \end{aligned}
\end{equation}

\begin{lm}[\cite{KM01}]\label{solRF}
{\it The Ricci flow on $E(1, 1)$ has the following properties.} 
% \begin{itemize}
\begin{enumerate}
\item {\it The solution of the Ricci flow equation $(\ref{sol3rf})$ exists on $(-T, +\infty)$, where $T>0$ 
is depending only on initial data.}
\item The quantity $a c$ and $b ( c-a )$ are conserved. 
\item {\it Put $\rho :=a/c$, then} 
\begin{equation*}
    \begin{aligned}
\dfrac{d}{dt} b & =\dfrac{(1+\rho )^2}{2 \rho } , \\
\dfrac{d}{dt} \rho & = \dfrac{1-\rho ^2}{b} .
	     \end{aligned}
\end{equation*}
\item $\rho _0 < \rho _{\infty} \le 1$ {\it or} $1\le \rho _{\infty} <\rho _0.$
\item {\it If $\rho \ne 1$, 
then $$b=k_0 \frac{\sqrt{\rho }}{|1-\rho |} , $$ where} $k_0= b_0\frac{|1-\rho _0|}{\sqrt{\rho _0}} .$ 
\item {\it If $t \rightarrow +\infty$, then $\rho \rightarrow 1$ and $b \rightarrow \infty$.}
\item {\it If $t \rightarrow +\infty$, then $a$ and $c$ converge to}
$a_{\infty}=c_{\infty}=\sqrt{a_0c_0}$.
% \end{itemize}
\end{enumerate}
\end{lm}
We use this proposition later to describe the asymptotic behavior of the resulting 
Ricci-flat metric on the space-time.

\begin{rem}
If $a = c$, then the metric $(\ref{sol})$ becomes a non-gradient expanding Ricci soliton. 
See Theorem $\ref{RSex}$.
\end{rem}

Next, let $(M^4, g)$ be a cohomogeneity one with respect to $E(1, 1)$.
Let $\{ F_i\} _{i=1}^3$ and $\{ \theta _i\} $ be as before.
 % Milnor frame of $\widetilde{{\rm Isom}({\bf E^1_1})}$,  \\
% $[F_2, F_3]=-2F_1, [F_3, F_1]=0, [F_1, F_2]=2F_3$. \\
% $\{ \theta _i\} :$ dual coframe of $\{ F_i\}$.
Then a cohomogeneity one metric is expressed as
\begin{equation}\label{sol_cohomo}
g = dt^2 + a (t)^2 (\theta ^1)^2 + b (t)^2 (\theta ^2)^2 + c (t)^2 (\theta ^3)^2 . 
\end{equation}
Putting $F_0 =\frac{\partial}{\partial t}$, we have an orthogonal frame $\{ F_i\}_{i=1}^{4}$. 
The Ricci tensor is computed as
\begin{subequations}
    \begin{align*}
  R_{00} & = -\frac{\ddot{a}}{a} -\frac{\ddot{b}}{b} -\frac{\ddot{c}}{c} , \\
R_{11} & = -a \frac{(\dot{a} b c)^{\cdot }}{b c} +\frac{a^4 - c^4}{2 b^2 c^2} , \\
R_{22} & = -b \frac{(\dot{b} c a)^{\cdot }}{c a} -\frac{(a^2 + c^2)^2}{2 a^2 c^2} , \\
R_{33} & = -c \frac{(\dot{c} a b)^{\cdot }}{a b} +\frac{c^4 - a^4}{2 a^2 b^2} ,
 \end{align*}
 \end{subequations}
and other component are $0$.
The Ricci flow equation for $(\ref{sol})$ is equivalent to the following: 
\begin{equation}\label{rfsol}
    \begin{aligned}
    \dfrac{d}{dt} a & =\dfrac{c^2 - a^2}{2 b c} , \\
    \dfrac{d}{dt} b & =\dfrac{(c + a)^2}{2 a c} , \\
    \dfrac{d}{dt} c & =\dfrac{a^2 - c^2}{2 a b} .
     \end{aligned}
\end{equation}
Since 
$ \dfrac{d}{dt} (a c) =0$,  it is easy to check that 
\begin{subequations}
    \begin{align*}
R_{22} & = -b \frac{(\dot{b} c a)^{\cdot }}{c a} -\frac{(c^2 + a^2)^2}{2 a^2 c^2} , \\
 & = -b \ddot{b} -\frac{(c^2 + a^2)^2}{2 a^2 c^2} , \\
 & =\frac{(c^2 - a^2)^2}{2 a^2 c^2} -\frac{(c^2 + a^2)^2}{2 a^2 c^2} = -2 , 
 \end{align*}
 \end{subequations}
and this metric is not a Ricci-flat metric.
%%%%ricci soliton, R.S. , Ricci soliton 
\begin{rem}
The metric $(\ref{sol_cohomo})$ with $(\ref{rfsol})$ is a Ricci soliton 
if and only if $a = -c$, in other words, $a$, $b$ and $c$ are constants.
Furthermore, this Ricci soliton is the product metric on $\Real \times E(1, 1)$ 
and a non-gradient expanding Ricci soliton.
We can check in the following way. 
First, because $E(1, 1) \approx \Real ^3$, 
we can use standard coordinates $(x_1, x_2, x_3)$ on $\Real ^3$.
Secondly, since left-invariant vector fields is given by 
$$
F_1 = e^{x_2} \frac{\partial }{\partial x_3} +e^{-x_2} \frac{\partial }{\partial x_1} , \ 
F_2 = \frac{\partial }{\partial x_2} , \ 
F_3 = e^{x_2} \frac{\partial }{\partial x_3} -e^{-x_2} \frac{\partial }{\partial x_1} , 
$$
then we can write the Ricci soliton equation using standard coordinates $(x_1, x_2, x_3)$  on $\Real ^3$, 
and we will obtain conditions of the Ricci soliton, for example, 
$$X = \frac{\alpha t}{2}  F_0 + e^{-x_2} \frac{\alpha x_3}{4} F_1 +\frac{\alpha }{4} F_2 \ , \  \alpha = \frac{4}{b^2} .$$
Lastly, we check that $\nabla _i X_j$ is not symmetric. 
Therefore this Ricci soliton is a non-gradient expanding Ricci soliton. 
\end{rem}

Next, we consider a cohomogeneity one metric  
\begin{equation}\label{sol22}
g=dt^2 + a(t)^2 (\theta ^1)^2 - b(t)^2 (\theta ^2)^2 - c(t)^2 (\theta ^3)^2. 
\end{equation}
This metric has signature $(2, 2)$. 
Put $F_0 =\frac{\partial}{\partial t}$. Then $\{ F_i\} $ is the orthogonal frame.
The Ricci tensor is given by
\begin{subequations}
    \begin{align*}
R_{00} & = -\frac{\ddot{a}}{a} -\frac{\ddot{b}}{b} -\frac{\ddot{c}}{c} , \\
R_{11} & = -a \frac{(\dot{a} b c)^{\cdot }}{b c} +\frac{a^4 - c^4}{2 b^2 c^2} , \\
R_{22} & = b \frac{(\dot{b} c a)^{\cdot }}{c a} +\frac{(c^2 - a^2)^2}{2 a^2 c^2} , \\
R_{33} & = c \frac{(\dot{c} a b)^{\cdot }}{a b} +\frac{a^4 - c^4}{2 a^2 b^2} .
 \end{align*}
 \end{subequations}
These formulas imply the following theorem: 
\begin{thm}
{\it Let $a$, $b$ and $c$ satisfy the Ricci flow equations; 
\begin{equation}
    \begin{aligned}
    \dfrac{d}{dt} a & =\dfrac{c^2 - a^2}{2 b c} , \\
    \dfrac{d}{dt} b & =\dfrac{(c + a)^2}{2 c a} , \\
    \dfrac{d}{dt} c & =\dfrac{a^2 - c^2}{2 a b} . 
     \end{aligned}
\end{equation}
Then the metric $(\ref{sol22})$ becomes a Ricci-flat metric.}
\end{thm}

\begin{proof}
\begin{equation}
    \begin{aligned}
    \dfrac{d^2}{dt^2} a & = \dfrac{(a^2 - c^2)(2 a^2+ a c+ c^2)}{a b^2 c^2} , \\
    \dfrac{d^2}{dt^2} b & = -\dfrac{(c^2-a^2)^2}{2 a^2 b c^2} , \\
    \dfrac{d^2}{dt^2} c & = \dfrac{(a^2- c^2)(a^2 + a c +2 c^2)}{2 a^2 b^2 c} ,
     \end{aligned}
\end{equation}
and $a c$ is a constant. 
Using this fact, we obtain the statement.
\end{proof}
%%%keisan is 2times.

Sectional curvatures of the metric $(\ref{sol22})$ satisfying $(\ref{sol3rf})$ are computed as
\begin{subequations}
    \begin{align*}
K_{01} =K_{23} & = \frac{(c^2-a^2)(2 a^2 +a c+c^2)}{2 a^2 b^2 c^2}  
= \frac{(1-\rho ^2)(2\rho ^2+\rho +1)}{2 \rho  ^2 b^2} , \\
K_{02} =K_{13} & = \frac{(c^2-a^2)^2}{2 a^2 b^2 c^2}  = \frac{(1-\rho ^2)^2}{2 \rho ^2 b^2} , \\
K_{03} =K_{12}& = -\frac{(c^2-a^2)(a^2+a c+2 c^2)}{2 a^2 b^2 c^2} 
= -\frac{(1-\rho ^2)(\rho ^2+\rho +2)}{2 \rho ^2 b^2} , 
   \end{align*}
  \end{subequations}
where $\rho := a/c $.
Notice that if $a^2 \ne c^2$, then signatures of sectional curvatures $K_{ij}$ are decided.
From Proposition $(\ref{solRF})$, we get the following proposition.

\begin{prop}
{\it Let $a$, $b$ and $c$ satisfy the Ricci flow equations $(\ref{sol3rf})$. 
% \begin{itemize}
\begin{enumerate}
\item If $a^2 = c^2$, then the metric $(\ref{sol22})$ becomes a zero curvature metric. 
\item If $a^2 \neq c^2$, then the metric $(\ref{sol22})$ is a Ricci-flat metric, but $K_{ij}$ are 
not constant.
\end{enumerate}
Assume that $a_0 c_0$ and $b_0$ are positive.
\begin{enumerate}
\item
The Ricci flow exists on $( -T, \infty )$, where $T$ depends on the initial data.
\item 
If $t \rightarrow +\infty$,  
then all sectional curvatures $K_{ij}  \rightarrow 0$.
\item
If $t \rightarrow -T$,  
then $K_{ij}   \rightarrow  \pm \infty $.
\end{enumerate}
Assume that $a_0c_0$ is positive, and $b_0$ is negative.
\begin{enumerate}
\item
The Ricci flow exists on $( -\infty ,T)$, where $T$ depends on initial data.
\item 
If $t \rightarrow -T$,  
then sectional curvature $K_{ij}  \rightarrow 0$.
\item
If $t \rightarrow \infty $, 
then $K_{ij} \rightarrow   \pm \infty $.
\end{enumerate}
Assume that $a_0 c_0$ is negative, and $b_0$ is positive.
Then the behavior of $a$, $b$ and $c$ is the backward Ricci flow of $E(2)$ 
(see $(\ref{e2})$). 

Assume that $a_0 c_0$ and $b_0$ are negative. 
Then the behavior of $a$, $b$ and $c$ is the Ricci flow of $E(2)$ (see $(\ref{e2})$).}
\end{prop}

% \begin{rem} %teisei
% For a metric of signature $(2, 2)$, a plane spanned by $X$ and $Y$ collapse, 
% if and only if $g(X, X)\cdot g(Y,Y) - g(X, Y)^2 = 0$, where $X$ and $Y$ are tangent vectors.
% A plane span $\{ X, Y\}$ is not collapse, if and only if
% $g(X, X)\cdot g(Y,Y) - g(X, Y)^2 \ne 0$.

% If a plane span $\{ X, Y\}$ is non-collapsing, 
% then we can define the sectional curvature of plane of $X$ and $Y$.
% Let a plane of $X$ and $Y$ be non-collapsing.
% Then the sectional curvature of the plane spanned $\{ X, Y\}$ is defined as
% $$K(X, Y) = \frac{R(X, Y, Y, X)}{g(X, X)\cdot g(Y,Y) - g(X, Y)^2} \ .$$
% This description depends only on the plane. 
% The fact is similar to the Riemann case.

% In this case, a collapsing plane is 
% $X=c F_0 + F_3, Y=F_2$, for example.
% \end{rem}

 %%%%%%%%%%%%%%%%%%%%$\widetilde{{\rm Isom}({\bf E^2})}$}-Geometry, sol
\section{The group of rigid motions of the Euclidean $2$-space}
In this section, we consider the group $E(2)$ of rigid motions of the Euclidean $2$-space.
Let $\{ F_i\}_{i=1}^3$ be the Milnor frame of $E(2)$, 
satisfying 
$$[F_2, F_3] = -F_1 \ , \ [F_3, F_1] = 0 \ , \ [F_1, F_2] = -F_3,$$
 and  $\{ \theta _i\}$ dual coframe of $\{ F_i\}$. 
Then a left-invariant metric is expressed as
$$g=a (\theta ^1)^2+b(\theta ^2)^2+c (\theta ^3)^2 . $$
Therefore the Ricci flow equation is equivalent to
\begin{equation}\label{e2}
    \begin{cases}
    \dfrac{d}{dt} a =\dfrac{c^2 - a^2}{2 b c} , \\
    \dfrac{d}{dt} b =\dfrac{(c - a)^2}{2 c a} , \\
    \dfrac{d}{dt} c =\dfrac{a^2 - c^2}{2 a b} .
     \end{cases}
\end{equation}
Since $E(1, 1)$ and $E(2)$ are solvable,
so the Ricci flow equations (\ref{e2}) resemble as the Ricci flow equations $(\ref{sol3rf})$.
The behavior of the solution of $(\ref{e2})$ is known in \cite{KM01, CK}. 
\begin{lm}[\cite{KM01, CK}]
{\it The Ricci flow on $E(2)$ has the following properties.}
\begin{enumerate}
\item {\it The solution of the the Ricci flow equation $(\ref{e2})$ exists on $(-T, \infty)$, 
where $T>0$ is depending only on initial data.}
\item The quantities $a c$ and $b ( c + a ) $ are preserved. 
\item {\it Put $k := a/c$. Then we get } $$\dfrac{d}{dt} b =\dfrac{(1-k )^2}{2 k} , $$
	$$\dfrac{d}{dt} k = \dfrac{1 - k^2}{b} .$$ 
\item {\it If $t \rightarrow \infty$, then $k \rightarrow 1$ and} $b \rightarrow ^{\exists} b_{\infty}$.
\item {\it The coefficient $b$ satisfies
$$b = l_0 \frac{\sqrt{k}}{ 1 + k } , $$ where} $ l_0 = b_0 \frac{1 + k_0}{\sqrt{k_0}}$ .
\end{enumerate}
\end{lm}
We use above proposition later to describe the asymtptotic property of the resulting 
Ricci-flat metric on the space-time. 
\begin{rem}
Since $E(2)$ is dynamically stable, 
the Ricci flow on $E(2)$ converges to an Einstein metric.
The details are written in \cite{IJ92, S06, SSS08}.
\end{rem}
%%%%kensyou
Next, we consider a cohomogeneity one metric with respect to $E(2)$. 
Let $\{ F_i\} _{i=1}^3$ and  $\{ \theta _i\}$ be as before.
% : Milnor frame of $\widetilde{{\rm Isom}({\bf E^2})}$,  
%  $[F_2, F_3]=-2F_1, [F_3, F_1]=0, [F_1, F_2]=-2F_3$. 
%  $\{ \theta _i\}$ : dual coframe of $\{ F_i\}$. 
Then the metric is expressed as
\begin{equation}\label{cohomo_e2}
g=dt^2+a(t)^2 (\theta ^1)^2+b(t)^2 (\theta ^2)^2+c(t)^2 (\theta ^3)^2 . 
\end{equation}
The Ricci tensor of this metric is given by
\begin{subequations}
    \begin{align*}
R_{00} & = -\frac{\ddot{a}}{a} -\frac{\ddot{b}}{b} -\frac{\ddot{c}}{c} , \\
R_{11} & = -a \frac{(\dot{a} b c)^{\cdot }}{b c} + \frac{a^4 - c^4}{2 b^2 c^2} , \\
R_{22} & = -b \frac{(\dot{b} c a)^{\cdot }}{c a} - \frac{(a^2 - c^2)^2}{2 a^2 c^2} , \\
R_{33} & = -c \frac{(\dot{c} a b)^{\cdot }}{a b} + \frac{c^4 - a^4}{2 a^2 b^2} ,
 \end{align*}
 \end{subequations}
 and other components are $0$.
Using the solution of $(\ref{e2})$, we construct a Ricci-flat metric on the space-time in the following way: 
\begin{thm}
{\it If $a$, $b$ and $c$ satisfy the Ricci flow equations 
\begin{equation}\label{rf_e2}
    \begin{cases}
    \dfrac{d}{dt} a =\dfrac{c^2 - a^2}{2 b c} , \\
    \dfrac{d}{dt} b =\dfrac{(a - c)^2}{2 a c} , \\
    \dfrac{d}{dt} c =\dfrac{a^2 - c^2}{2 a b} , 
     \end{cases}
\end{equation}
then this metric $(\ref{cohomo_e2})$ becomes a Ricci-flat metric. } 
\end{thm}

\begin{proof}
This is proved by routine calculations. 
\end{proof}
%%%%%%%%%%%%%%%keisan is 2 times.
Sectional curvatures of this metric are
\begin{subequations}
    \begin{align*}
K_{01} =K_{23} & = \frac{(c^2 - a^2)(2 a^2-a c + c^2)}{a^2 b^2 c^2} = \frac{(1-k^2) (2 k^2- k+1)}{2 b^2 k^2} , \\
K_{02} =K_{13} & = \frac{(c^2 - a^2)^2}{2 a^2 b^2 c^2} = \frac{(1 -k^2)^2}{2 b^2 k^2} , \\
K_{03} =K_{12} & = -\frac{(c^2 - a^2)(a^2- c a + 2 c^2)}{2 a^2 b^2 c^2} = - \frac{(1 - k^2)(k^2- k+2)}{2 b^2 k^2} , 
    \end{align*}
  \end{subequations}
where $k := a/c$.
We can observe similarities between the Ricci flow solutions on the space-times of the 
left-invariant Ricci flow solutions of $E(1, 1)$ and $E(2)$. 
For instance, we get the following proposition. 
\begin{prop}
\begin{enumerate} 
\item {\it If $a^2 = c^2$, then this metric is a flat curvature metric. } 
\item {\it If $a^2 \neq c^2$, then this metric is a Ricci-flat metric but not a constant curvature metric.} 
\end{enumerate}
{\it Let $a_0c_0$ and $b_0$ be positive.}
\begin{enumerate} 
\item {\it The Ricci flow equation exists on $( -T, \infty )$, where $T$ depends on initial data.}
\item {\it If} $t\rightarrow +\infty$, {\it then} $K_{ij}\rightarrow 0$. 
\end{enumerate}
{\it Let $a_0 c_0$ be positive, and $b_0$ be negative.}
\begin{enumerate}
\item
{\it The Ricci flow equation exists on $( -\infty ,T)$, where $T$ depends on initial data.}
\item 
{\it If $t\rightarrow -T$,  
then sectional curvatures $K_{ij} \rightarrow 0$.}
\item
{\it If $t\rightarrow \infty $, 
then $K_{ij} \rightarrow  \pm \infty $.}
\end{enumerate}
{\it Let $a_0 c_0$ be negative, and $b_0$ be positive.
Then behavior of $a$, $b$ and $c$ is the backward Ricci flow of $E(1, 1)$ (see $(\ref{sol3rf})$).
Let $a_0 c_0$ and $b_0$ be negative. 
Then behavior of $a$, $b$ and $c$ is the Ricci flow of $E(1, 1)$ (see $(\ref{sol3rf})$). }
\end{prop}

Next,  we define almost complex structures.
Using $E(2) \approx \Real ^3$, we can write
\begin{subequations}
    \begin{align*}
F_1 & = \sin y \frac{\partial }{\partial x} +\cos y \frac{\partial }{\partial z} , \\
F_2 & = \frac{\partial }{\partial y} , \\
F_3 & = \cos y \frac{\partial }{\partial x} -\sin y \frac{\partial }{\partial z} ,
 \end{align*}
 \end{subequations}
and 
\begin{subequations}
    \begin{align*}
\theta ^1 & = \sin y \cdot dx +\cos y \cdot dz , \\
\theta ^2 & = dy , \\
\theta ^3 & = \cos y \cdot dx -\sin y \cdot dz .
 \end{align*}
 \end{subequations}
Then the metric on $E(2)$ with $a = b = c = 1$ is
$$g_3 = (\theta ^1)^2 + (\theta ^2)^2+ (\theta ^3)^2 = d x^2 + d y^2 + d z^2 . $$
Let $\{ e_i\} _{i = 0}^3$ be an orthonormal coframe defined by 
$$e^0 = d t \ , \  e^1=a\theta ^1 \ , \  e^2=b\theta ^2 \ , \ e^3=c\theta ^3 .$$
Then $\{ e^i\} _{i=0}^4$ becomes an orthonormal coframe of cohomogeneity one metric $g$.
Let $\{ e_i\}$ be the orthonormal frame, satisfying $e^i (e_j) =\delta ^i_j$.
Since $\Real ^4 \approx \Complex ^2$ has natural almost complex structures $\{ J_i\}$,
we analogously define almost complex structures $J_1$ and $J_2$, satisfying
\begin{subequations} %teisei
    \begin{align*}
J_1 e_0 & =e_2,  \quad J_1 e_3 =e_1 , \quad (J_1) ^2=-id , \\
J_2 e_0 & =\cos y \cdot e_1 -\sin y \cdot e_3 , \\
J_2 e_2 & =\sin y \cdot e_1 + \cos y \cdot e_3 , \\
(J_2) ^2  & =-id .
\end{align*}
 \end{subequations} %teisei
If we put $J_3:=J_2 J_1 =-J_1 J_2$, then $J_3$ becomes an almost complex structure. 
From direct computation we have $N_{J_i}(F_j,F_k)=0$ for all $i,j,k=0,1,2,3$ where 
$N_{J_i}$ is the Nijenhuis tensor of $J_i$. Therefore $J_i$'s are all integrable. 
We consider a triple of 2-forms $(\omega ^1, \omega ^2, \omega ^3)$, defined by 
$$\omega ^i (X, Y) := g(J_i X, Y)$$
for $i = 1, 2, 3$.
Then 2-forms $\omega _i$ are given by
\begin{subequations}
    \begin{align*}
\omega ^1 & =e^0 \wedge e^2 + e^3 \wedge e^1 , \\
\omega ^2 & =\cos y \cdot (e^0 \wedge e^1 + e^2 \wedge e^3) -\sin y \cdot (e^0 \wedge e^3 + e^1 \wedge e^2) , \\
\omega ^3 & =\sin y \cdot (e^0 \wedge e^1 + e^2 \wedge e^3) +\cos y \cdot (e^0 \wedge e^3 + e^1 \wedge e^2) .
\end{align*}
 \end{subequations}
 The above triple of $2$-forms $(\omega ^1, \omega ^2, \omega ^3)$ satisfies the following relations: 
$$(\omega ^i)^2 \ne 0  , \  \omega ^i \wedge \omega ^j =0 \ (i \ne j) .$$
Since $d e^0 =d ( d t ) = 0 , \  d e^1 = \dot{a} d t \wedge \theta ^1 + a\theta ^2 \wedge \theta ^3 ,
 \ de^2 = \dot{b} dt \wedge \theta ^2 , \ d e^3=\dot{c} dt \wedge \theta ^3 + c\theta ^1 \wedge \theta ^2 ,  $
 we get 
\begin{subequations} 
    \begin{align*}
d\omega ^1 & = (\dot{c}a +c\dot{a} ) dt\wedge \theta ^3 \wedge \theta ^1 , \\
d\omega ^2 & =-(\dot{a}b +a\dot{b} + a- c )\sin y\cdot dt\wedge \theta ^1 \wedge \theta ^2
		+(\dot{b} c + b \dot{c} + c- a )\cos y\cdot dt\wedge \theta ^2 \wedge \theta ^3 , \\
d\omega ^3 & = - (\dot{a}b +a\dot{b} + a- c )\cos y\cdot dt\wedge \theta ^1 \wedge \theta ^2
		+(\dot{b} c + b \dot{c} + c- a )\sin y\cdot dt\wedge \theta ^2 \wedge \theta ^3 . 
\end{align*}
 \end{subequations}

\begin{theorem}
{\it The cohomogeneity one metric $g$ expressed as $(\ref{cohomo_e2})$, where the triple 
$\{a(t),b(t),c(t)\}$ satisfies the Ricci flow equation on $E(2)$, 
is a hyper-K\"{a}hler metric. 
Conversely the cohomogeneity one metric $g$ as in $(\ref{cohomo_e2})$ 
satisfies  $d\omega ^i =0$ for any $i=1, 2, 3$, 
then the triple $\{a(t),b(t),c(t)\}$ satisfies the Ricci flow equations on $E(2)$ . }
\end{theorem}

\proof 
It is straitforward to check that 
$d\omega ^i =0$ if and only if the following hold: 
\begin{subequations}
    \begin{align}
\dot{c} a+ c\dot{a} & = 0  \label{ia} , \\
\dot{a} b +a\dot{b}  & = c - a \label{ic} , \\
\dot{b} c +b\dot{c}  & = a - c \label{ib} .
\end{align}
\end{subequations}
$(\ref{ia} )$ implies that $c a$ is constant. 
Solveing the equations for $(\dot{a}, \dot{b} ,\dot{c})$, 
we obtain the statement of the theorem. 
\endproof
Therefore the cohomogeneity one metric $(\ref{rf_e2} )$ becomes a hyper-K\"{a}hler metric.
\begin{rem} %teisei 
In this section, we obtain the cohomogeneity one metric with respect to $E(2)$ with the hyper-K\"{a}hler structure $\{ J_i\}$.
From the definition of $\{ J_i\}$, the hyper-K\"{a}hler structure is non-trivial. 
The cohomogeneity one metric with respect to $E(2)$ with the trivial hyper-K\"{a}hler structure (i.e. $J_i e_0 = e_i, J_i e_j = e_k$, where $(i, j, k)$ are cyclic permutation of $\{1,2,3\}$,)
is also a hyper-K\"{a}hler structure, 
however the triple of functions $\{a(t),b(t),c(t)\}$ does not satisfy the Ricci flow equation.  
So we do not consider it deeply.  %teisei 
\end{rem}

%%%%%teisei, tekito
% \begin{rem}
% The metric $(\ref{cohomo_e2})$ with $(\ref{rf_e2})$ are linear stable iff $a^2=c^2$.
% If $a^2 \ne c^2$, then the metric $(\ref{cohomo_e2})$ with $(\ref{rf_e2})$ is $($ not $?)$ linear stable.
% \begin{proof}
% Suppose $a^2 >c^2$. 
% Then signature of sectional curvature are as follows, 
% $K_{01}, K_{13}, K_{23} >0$ and $K_{02}, K_{03}, K_{12} <0$.
% \end{proof}
% \end{rem}
%%%%%%%%%%%%%%%%%%%%%%%%%%%%%%%%%%%$\widetilde{{\rm SL}(2, \Real )}$ sl

\section{$SL(2, \Real )$}
In this section, we consider $SL(2, \Real )$.
Let $(M^4, g)$ be a cohomogeneity one with respect to $SL(2, \Real )$. 
Let $\{ F_i\} _{i=1}^3$ be the Milnor frame of $SL(2, \Real )$, 
satisfying 
$$[F_1, F_2]= - F_3 , \  [F_2, F_3]= F_1 , \ [F_3, F_1] = F_2 , $$ 
and $\{ \theta _i\}$ the dual coframe of $\{ F_i\}$.
The metric is expressed as
\begin{equation}\label{cohomo_sl2}
g=dt^2+a(t)^2 (\theta ^1)^2+b(t)^2 (\theta ^2)^2+c(t)^2 (\theta ^3)^2 .
\end{equation}
The Ricci tensor of this metric is given by
\begin{subequations}
    \begin{align*}
  R_{00} & = -\frac{\ddot{a}}{a} -\frac{\ddot{b}}{b} -\frac{\ddot{c}}{c} , \\
R_{11} & = -a \frac{(\dot{a} b c)^{\cdot }}{b c} - \frac{(b^2 + c^2)^2-a^4}{2 b^2 c^2} , \\
R_{22} & = -b \frac{(\dot{b} c a)^{\cdot }}{c a} - \frac{(a^2 + c^2)^2-b^4}{2 a^2 c^2} , \\
R_{33} & = -c \frac{(\dot{c} a b)^{\cdot }}{a b} - \frac{(a^2 - b^2)^2-c^4}{2 a^2 b^2} ,
 \end{align*}
 \end{subequations}
and other component are $0$.

\begin{prop}
{\it If $a$, $b$ and $c$ satisfy the Ricci flow equations of $SL(2, \Real )$ 
\begin{equation}\label{rf_sl2}
    \begin{cases}
    \dfrac{d}{dt} a =\dfrac{(b +c)^2 -a^2}{2 b c} , \\
    \dfrac{d}{dt} b =\dfrac{(c +a)^2 -b^2}{2 c a} , \\
    \dfrac{d}{dt} c =\dfrac{(a -b)^2 -c^2}{2 a b} , 
     \end{cases}
\end{equation}
then this metric $(\ref{cohomo_sl2} )$ satisfying $(\ref{rf_sl2})$ is not a Ricci-flat metric.} 
\end{prop}

\begin{proof}
The proof is a straightforward calculation. 
\end{proof}

Next, we consider the metric of the form
\begin{equation}\label{sl2_22} 
g=dt^2 - a (t)^2 (\theta ^1)^2 - b (t)^2 (\theta ^2)^2 + c (t)^2 (\theta ^3)^2. 
\end{equation}
This metric $(\ref{sl2_22})$ has signature $(2, 2)$.

% \begin{eqnarray*}
% R(F_0, F_1, F_1, F_0) & = & a\ddot{a} , \\
% R(F_0, F_2, F_2, F_0) & = & b\ddot{b} , \\
% R(F_0, F_3, F_3, F_0) & = & -c\ddot{c} , \\
% R(F_1, F_2,F_2, F_1) & = & -a\dot{a} b\dot{b} +\frac{(a^2-b^2)^2-c^4}{c^2} -2(a^2+b^2+c^2), \\
% R(F_1, F_3,F_3, F_1) & = & c\dot{c} a\dot{a} -\frac{(a^2+c^2)^2-b^4}{b^2} +2(a^2-b^2-c^2), \\
% R(F_2, F_3,F_3, F_2) & = & b\dot{b} c\dot{c} -\frac{(b^2+c^2)^2-a^4}{a^2} -2(a^2-b^2+c^2), 
%  \end{eqnarray*} 

The Ricci tensor of this metric is computed as
\begin{subequations}
    \begin{align*}
R_{00} & = -\frac{\ddot{a}}{a} -\frac{\ddot{b}}{b} -\frac{\ddot{c}}{c} , \\
R_{11} & = a \frac{(\dot{a} b c)^{\cdot }}{b c} + \frac{(b^2 - c^2)^2-a^4}{2 b^2 c^2} , \\
R_{22} & = b \frac{(\dot{b} c a)^{\cdot }}{c a} + \frac{(a^2 - c^2)^2-b^4}{2 a^2 c^2} , \\
R_{33} & = -c \frac{(\dot{c} a b)^{\cdot }}{a b} - \frac{(a^2 - b^2)^2-c^4}{2 a^2 b^2} ,
  \end{align*}
 \end{subequations}
 and other components are $0$.
The Ricci flow equations of $SL(2, \Real )$ are given by
\begin{equation} 
    \begin{aligned}
    \dfrac{d}{dt} a =\dfrac{(b + c)^2-a^2}{2 b c} , \\
    \dfrac{d}{dt} b =\dfrac{(c + a)^2-b^2}{2 c a} , \\
    \dfrac{d}{dt} c =\dfrac{(a - b)^2-c^2}{2 a b} .
     \end{aligned}
\end{equation}
In this case,  
changing $a$ into $-a$ and $b$ into $-b$,
therefore the Ricci flow equations of $SL(2, \Real )$ change into
\begin{equation}\label{su2rf}
    \begin{aligned}
    \dfrac{d}{dt} a =\dfrac{(b - c)^2 - a^2}{2 b c} , \\
    \dfrac{d}{dt} b =\dfrac{(c - a)^2 - b^2}{2 c a} , \\
    \dfrac{d}{dt} c =\dfrac{(a - b)^2 - c^2}{2 a b} .
     \end{aligned}
\end{equation}
This system of equations is nothing but the Ricci flow equations of left-invariant metrics 
on $SU(2)$.

\begin{thm}
If  $a$, $b$ and $c$ satisfy the Ricci flow equations of $SU(2)$, 
then the metric $(\ref{sl2_22})$ satisfying $(\ref{su2rf})$ becomes a Ricci-flat metric 
of signature $(2,2)$.
\end{thm}

\begin{proof}
The proof method of this theorem is similar to Theorem $\ref{su2ricciflat}$.
The proof is straightforward. 
\end{proof}

\end{document}